\documentclass[a4paper,11pt]{article}
\usepackage{latexsym}
\usepackage{amsmath}
\usepackage{amssymb}
\usepackage{enumerate}
\usepackage{theorem}
\usepackage{array}
\newtheorem{itheorem}{Theorem}

\newtheorem{theorem}{Theorem}[section]
\newtheorem{lemma}[theorem]{Lemma}
\newtheorem{corollary}[theorem]{Corollary}
\newtheorem{proposition}[theorem]{Proposition}
\theorembodyfont{\rmfamily}
\newtheorem{definition}[theorem]{Definition}

\newtheorem{setting}[theorem]{Setting}

\newtheorem{notationinintro}{Notation}

\newcommand{\proof}{\noindent \mbox{\em Proof.\hspace*{2mm}}}
\newcommand{\qed}{\hfill \mbox{$  \Box $}}
\DeclareMathOperator{\Sing}{Sing}
\DeclareMathOperator{\image}{Im}
\newcommand{\gyokan}{\vskip 6pt}
\title{Numbers of points of surfaces in the projective $3$-space
over finite fields
\footnote{submitted to Finite Fields and their Applications.}}
\author{
Masaaki Homma
\thanks{Partially supported by Grant-in-Aid
for Scientific Research (24540056), JSPS.}
\\
 Department of Mathematics and Physics\\
Kanagawa University\\
Hiratsuka 259-1293, Japan\\
homma@kanagawa-u.ac.jp
\and
Seon Jeong Kim
\thanks{Partially supported by Basic Science Research Program through the National Research Foundation of Korea(NRF) 
funded by the Ministry of Education, Science and Technology (2012R1A1A2042228),
and also by the Gyeongsang National University Fund for Professors on Sabbatical Leave, 2013.
}\\
 Department of Mathematics and RINS\\
Gyeongsang National University\\
Jinju 660-701, Korea \\
skim@gnu.kr
}
\date{}

\begin{document}
\maketitle
\begin{abstract}
In the previous paper,
we established an elementary bound for numbers of points of surfaces in
the projective $3$-space over ${\Bbb F}_q$.
In this paper, we give the complete list of surfaces that attain the elementary bound.
Precisely those surfaces are the hyperbolic surface, the nonsingular Hermitian surface, and the surface of minimum degree containing all ${\Bbb F}_q$-points of the $3$-space.
\end{abstract}

\section{Introduction}
Let $S$ be a surface of degree $d$ in ${\Bbb P}^3$ over ${\Bbb F}_q$
without ${\Bbb F}_q$-plane components,
and $N_q(S)$ the number of ${\Bbb F}_q$-points of $S$.
In the previous paper \cite{hom-kim2013b},
we established the elementary bound for $N_q(S)$:
\begin{equation}\label{elementarybound}
 N_q(S) \leq (d-1)q^2 + dq +1,
\end{equation}
and also gave three examples of surfaces that achieve the upper bound
(\ref{elementarybound}).
The goal of this paper is to show that
only those three are examples of such surfaces.
\begin{itheorem}
For a surface $S$ in ${\Bbb P}^3$ over ${\Bbb F}_q$ without ${\Bbb F}_q$-plane
components,
if equality holds in {\rm (\ref{elementarybound})},
then the degree $d$ of $S$ is either
$2$ or $\sqrt{q} +1$ {\rm(}when $q$ is a square{\rm)} or $q+1$.
Furthermore,
the surface $S$ is projectively equivalent to
one of the following surfaces over ${\Bbb F}_q${\rm :}
\begin{enumerate}[{\rm (i)}]
\item
$X_0X_1-X_2X_3=0$  if $d=2${\rm ;}
\item
$X_0^{\sqrt{q}+1}+ X_1^{\sqrt{q}+1} +X_2^{\sqrt{q}+1}+X_3^{\sqrt{q}+1}=0$
if $d= \sqrt{q}+1${\rm ;}
\item $X_0X_1^{q} - X_0^{q}X_1 + X_2X_3^{q} - X_2^{q}X_3 =0$
if $d=q+1$.
\end{enumerate}
\end{itheorem}
\begin{notationinintro}
For an algebraic set $X$ defined by equations over ${\Bbb F}_q$
in a projective space,
the set of ${\Bbb F}_q$-points of $X$ is denoted by $X({\Bbb F}_q)$,
and the cardinality of $X({\Bbb F}_q)$ by $N_q(X)$.
The symbol $\theta_q(r)$ denotes $N_q({\Bbb P}^r)$,
and we understand $\theta_q(0)=1$.

The set of ${\Bbb F}_q$-planes of ${\Bbb P}^3$ is denoted by
$\check{\Bbb P}^3({\Bbb F}_q)$.
For an ${\Bbb F}_q$-line $l$ in ${\Bbb P}^3$,
$\check{l}({\Bbb F}_q)$
denotes the set
$\{ H \in \check{\Bbb P}^3({\Bbb F}_q) \mid H \supset l \}.$

When $Y$ is a finite set, ${}^{\#}Y$ denotes the cardinality of $Y$.

When $M$ is a matrix, ${}^t M$ denotes the transposed matrix of $M$.
\end{notationinintro}
\section{Review of some results in our previous works}
\subsection{Plane curves}
To investigate plane sections of $S$, we need some results on plane curves.
\begin{proposition}[Sziklai bound]
Let $C$ be a curve of degree $d$ in ${\Bbb P}^2$ over ${\Bbb F}_q$
without ${\Bbb F}_q$-line components.
Then
\begin{equation}\label{sziklaibound}
 N_q(C) \leq (d-1)q +1
\end{equation}
unless $C$ is the curve over ${\Bbb F}_4$ defined by
\begin{equation}\label{exceptionalcase}
(X+Y+Z)^4+(XY+YZ+ZX)^2 +XYZ(X+Y+Z)=0
\end{equation}
with a certain choice of coordinates $X, Y, Z$ of ${\Bbb P}^2$.
For the exceptional case, the number of ${\Bbb F}_4$-points is $14$.
\end{proposition}
\proof
The proof is spreaded over three papers
\cite{hom-kim2009, hom-kim2010a, hom-kim2010b}.

\gyokan

Note that the bound (\ref{sziklaibound}) makes sense only for $d \leq q+2$
because $(d-1)q +1 > \theta_q(2)$ if $d >q+2$.

\begin{lemma}\label{7bitangents}
The curve {\rm (\ref{exceptionalcase})} has
$7$ bitangent lines, and those lines are defined over ${\Bbb F}_2$,
in particular they are also ${\Bbb F}_4$-lines.
\end{lemma}
\proof See \cite[Remark 3.1]{hom-kim2009}.

\gyokan

Although the classification of curves that attain the Sziklai bound
is still under way \cite{hom-kim2011, hom-kim2012},
several properties of such curves are known.
In the next lemma,
$\Sing C$ denotes the set of singular points of $C$.
\begin{lemma}\label{mustbeirreducible}
Let $C$ be a plane curve of degree $d$ over ${\Bbb F}_q$
without ${\Bbb F}_q$-line components such that
$N_q(C) = (d-1)q +1$.
Then $C$ is absolutely irreducible
and $\Sing C \cap C({\Bbb F}_q) = \emptyset$.
Furthermore, $d=2$ or $\sqrt{q} +1 \leq d \leq q+2$.
\end{lemma}
\proof
For the first part of the assertion, see \cite[\S 2]{hom-kim2010a}.
The Hasse-Weil bound for a plane curve which may have singularities
holds in the following form
\cite[Corollary 2.5]{aub-per1993}, \cite[Corollary 2]{lee-yeo1994}:
for an absolutely irreducible plane curve $C$ of degree $d$
over ${\Bbb F}_q$,
\begin{equation}\label{hasseweil}
 N_q(C) \leq q+1 + (d-1)(d-2)\sqrt{q}.
\end{equation}
The bound (\ref{hasseweil}) is better than (\ref{sziklaibound})
if $2 < d < \sqrt{q}+1$,
in fact,
$(d-1)q +1 -
\left(
q+1 + (d-1)(d-2)\sqrt{q}
\right)
= (d-2)(\sqrt{q} + 1 -d)\sqrt{q}.$
\qed
\subsection{Space surfaces}
Since the right-hand side of the inequality (\ref{elementarybound})
is bigger than $N_q({\Bbb P}^3) = \theta_q(3)$
if $d > q+1$,
the meaningful range for $d$ is $2 \leq d \leq q+1$.
From now on, we keep the following situation.
\begin{setting}\label{oursetting}
Let $d$ be an integer with $2 \leq d \leq q+1$,
and $S$ a surface of degree $d$ in ${\Bbb P}^3$
defined over ${\Bbb F}_q$
without ${\Bbb F}_q$-plane components.
We assume that
$N_q(S) = (d-1)q^2 + dq +1.$
\end{setting}

Under this setting,
we already observed several properties of lines on $S$
in \cite[\S 3]{hom-kim2013preprint}, especially
the surface in Setting~\ref{oursetting}
has an ${\Bbb F}_q$-line \cite[Lemma 3.3]{hom-kim2013preprint}.

\begin{definition}
Let $l_1, \dots , l_d$ be ${\Bbb F}_q$-lines in ${\Bbb P}^3$
with $d \geq 2$, and $Z = \cup_{i=1}^{d} l_i$.
The union of lines $Z$ is called
a planar ${\Bbb F}_q$-pencil of degree $d$
if those $d$ lines lie on a plane simultaneously
and they meet together at a point,
which is called the vertex of $Z$ and denoted by $v_Z$.

Since $d \geq 2$, the plane on which $Z$ lies and the vertex
of $Z$ are defined over ${\Bbb F}_q$.

For the surface $S$ in Setting~\ref{oursetting}, let
\[
\mathcal{Z}_S =
\{ H \in \check{\Bbb P}^3({\Bbb F}_q) \mid 
S\cap H \ \mbox{\rm  is a planar ${\Bbb F}_q$-pencil of degree $d$} \}.
\]
\end{definition}

\begin{lemma}\label{knownlemma1}
Let $l$ be an ${\Bbb F}_q$-line on the surface $S$
in Setting~{\rm \ref{oursetting}}.
\begin{enumerate}[{\rm (i)}]
\item If an ${\Bbb F}_q$-plane $H$ contains $l$, then
$S \cap H$ is a planar ${\Bbb F}_q$-pencil of degree $d$.
\item The map
$\check{l}({\Bbb F}_q) \ni H \mapsto v_{S\cap H} \in l({\Bbb F}_q)$
is bijective.
\end{enumerate}
\end{lemma}
\proof
See \cite[Lemma 3.6]{hom-kim2013preprint}.

\begin{lemma}\label{knownlemma2}
Let $S$ be the surface in Setting~{\rm \ref{oursetting}},
and $H$ an ${\Bbb F}_q$-plane such that $S\cap H$ is
a planar ${\Bbb F}_q$-pencil of degree $d$.
If an ${\Bbb F}_q$-line $l$ on $S$ goes through
the vertex $v_{S\cap H}$,
then $l$ lies on $H$, and so it is a component of $S\cap H$.
\end{lemma}
\proof
See \cite[Lemma 3.7]{hom-kim2013preprint}.

\begin{lemma}\label{knownlemma3}
Let $S$ be the surface in Setting~{\rm \ref{oursetting}},
and $P \in S({\Bbb F}_q)$.
Then there is an ${\Bbb F}_q$-plane $H$ containing $P$
such that $S\cap H$ is a planar ${\Bbb F}_q$-pencil
with $v_{S\cap H} = P$.
\end{lemma}
\proof
See \cite[Corollary 3.8]{hom-kim2013preprint}.

\begin{corollary}\label{knowncorollary}
Let $S$ be the surface in Setting~{\rm \ref{oursetting}}.
There is a natural bijection between
$S({\Bbb F}_q)$ and the set $\mathcal{Z}_S$ by
\[
\mathcal{Z}_S \ni H \mapsto v_{S\cap H} \in S({\Bbb F}_q).
\]
\end{corollary}
\proof
By Lemma~\ref{knownlemma3},
this map is surjective.
The injectivity comes from Lemma~\ref{knownlemma2}.
\qed

\section{Possible degrees}
The next proposition is a generalization of
\cite[Proposition 4.1]{hom-kim2013preprint} to any degree $d$.
\begin{proposition}\label{twocases}
Under Setting~{\rm \ref{oursetting}},
let $H$ be an ${\Bbb F}_q$-plane of ${\Bbb P}^3$.
Then $S \cap H$ is either
\begin{enumerate}[{\rm (1)}]
\item a planar ${\Bbb F}_q$-pencil of degree $d$, or
\item a plane curve of degree $d$ without ${\Bbb F}_q$-line
components with $N_q(S\cap H) = (d-1)q +1$.
\end{enumerate}
Furthermore, for $i= 1$ and $2$,
let $\nu_i$ denote the number of ${\Bbb F}_q$-planes
such that their sections on $S$ have the property $(i)$ above.
Then
$\nu_1= N_q(S)$ and $\nu_2 = \theta_q(3) - N_q(S)$.
\end{proposition}
\proof
If $S\cap H$ contains an ${\Bbb F}_q$-line,
then it is a planar ${\Bbb F}_q$-pencil of degree $d$ by
Lemma~\ref{knownlemma1}.
When $S\cap H$ does not contain any ${\Bbb F}_q$-line,
$N_q(S \cap H) \leq (d-1)q +1$ by the Sziklai bound,
except for the case $d=q=4$.
But even this case,
Sziklai's inequality holds under Setting~\ref{oursetting}.
In fact, suppose $d=q=4$ and $S \cap H_0$
were defined by (\ref{exceptionalcase}).
Then counting the number of ${\Bbb F}_q$-points of $S$ by using
a bitangent $l$ to $S \cap H_0$ in Lemma~\ref{7bitangents},
we would have
\begin{align*}
N_q(S) &= \sum_{H \in \check{l}({\Bbb F}_4)}
  \left( {}^{\#} ( S \cap H ({\Bbb F}_4)) -2 \right) +2 \\
      & \leq 5\cdot(14-2) +2 = 62,
\end{align*}
but $N_q(S)$ should be $(4-1)\cdot 4^2 + 4 \cdot 4 +1 =65.$

Consider the correspondence
\[
\mathcal{A} =
\{
(P, H) \in S({\Bbb F}_q) \times \check{\Bbb P}^3({\Bbb F}_q)
    \mid P \in H
\}
\]
together with projections
$\pi_1 : \mathcal{A} \to  S({\Bbb F}_q)$
and 
$\pi_2 : \mathcal{A} \to  \check{\Bbb P}^3({\Bbb F}_q)$.
Counting ${}^{\#}\mathcal{A}$ by using $\pi_1$,
we have
\begin{equation}\label{countingAfirst}
{}^{\#}\mathcal{A} = N_q(S)\, \theta_q(2).
\end{equation}
Recall
\begin{align*}
\mathcal{Z}_S &=
\{ H \in \check{\Bbb P}^3({\Bbb F}_q) \mid 
S\cap H \ \mbox{\rm  is a planar ${\Bbb F}_q$-pencil} \}\\
 &=
\{ H \in \check{\Bbb P}^3({\Bbb F}_q) \mid 
S\cap H \ \mbox{\rm  contains an  ${\Bbb F}_q$-line} \}.
\end{align*}
Then ${}^{\#} \mathcal{Z}_S = N_q(S)$ by Corollary~\ref{knowncorollary}.
Hence counting ${}^{\#}\mathcal{A}$
by using $\pi_2$,
we have
\begin{align}
{}^{\#}\mathcal{A} &= \sum_{H \in \mathcal{Z}_S} \pi_2^{-1}(H)
      + \sum_{H \in \check{\Bbb P}^3({\Bbb F}_q) \setminus \mathcal{Z}_S} 
                  \pi_2^{-1}(H) \notag \\
 & \leq N_q(S) (dq+ 1) + (\theta_q(3) - N_q(S))( (d-1)q +1) \tag{$*$}\\
 & =  N_q(S)q + \theta_q(3)( (d-1)q +1) \notag \\
 & = \left(  (d-1)q^2 + dq + 1  \right)q + (q^2+1)(q+1)( (d-1)q +1) \notag \\
 & = N_q(S)\, \theta_q(2). \notag
\end{align}
Taking account of (\ref{countingAfirst}),
we know that equality holds in ($*$),
which means that $N_q(S\cap H) = (d-1)q +1$ for any
$H \in {\Bbb P}^3({\Bbb F}_q) \setminus \mathcal{Z}_S.$
\qed

\begin{proposition}\label{possiblenumbers}
Under Setting~{\rm \ref{oursetting}},
let $l$ be an ${\Bbb F}_q$-line in ${\Bbb P}^3$.
Then ${}^{\#}(l \cap S({\Bbb F}_q))$ is either
$0$ or $1$ or $d$ or $q+1$.
\end{proposition}
\proof
For a given ${\Bbb F}_q$-line $l$,
the set $\{ H \in \mathcal{Z}_S \mid H \supset l \}$
is denoted by $\mathcal{Z}(l)$,
that is, $\mathcal{Z}(l) = \mathcal{Z}_S \cap \check{l}({\Bbb F}_q).$
Then ${}^{\#}(l \cap S({\Bbb F}_q))= {}^{\#}\mathcal{Z}(l)$.
To see this claim, put $\alpha = {}^{\#}(l \cap S({\Bbb F}_q))$
and $\beta = {}^{\#}\mathcal{Z}(l).$
Then
\begin{align*}
 N_q(S) &= \sum_{H \in \check{l}({\Bbb F}_q)}
                 (N_q(S \cap H) - \alpha) + \alpha \\
       &= \sum_{H \in \mathcal{Z}(l)}
             (N_q(S \cap H) - \alpha) +
          \sum_{H \in \check{l}({\Bbb F}_q) \setminus \mathcal{Z}(l)}
             (N_q(S \cap H) - \alpha) + \alpha \\
        & = \beta (dq+1 - \alpha) +(q+1-\beta)((d-1)q+1 -\alpha)
                   + \alpha \\
       & =\beta q +(d-1)q^2 + dq +1 - \alpha q.
\end{align*}
Since $N_q(S) = (d-1)q^2 + dq +1$,
we have $\alpha = \beta$.

Hence if ${}^{\#}(l \cap S({\Bbb F}_q)) >0$, there is an ${\Bbb F}_q$-plane
$H \in \check{l}({\Bbb F}_q)$ such that $S \cap H$ is a planar
${\Bbb F}_q$-pencil of degree $d$.
A planar ${\Bbb F}_q$-pencil of degree $d$ meets an ${\Bbb F}_q$-line
on the plane in $1$ or $d$ or $q+1$ ${\Bbb F}_q$-points.
Since $l \cap S = l \cap (S \cap H)$,
this completes the proof.
\qed

\begin{theorem}
In Setting~{\rm \ref{oursetting}},
the degree $d$ of $S$ must be either $q+1$ or $\sqrt{q}+1$ or $2.$
\end{theorem}
\proof
Suppose $d<q+1$.
Then $\theta_q(3) - N_q(S) >0$.
Hence there is an ${\Bbb F}_q$-plane $H$ such that
$S \cap H$ has no ${\Bbb F}_q$-line components and
$N_q(S \cap H) = (d-1)q +1$ by Proposition~\ref{twocases}.
Put $C=S \cap H$.
By Lemma~\ref{mustbeirreducible},
$C$ is absolutely irreducible and
$\Sing C \cap C({\Bbb F}_q) = \emptyset$.

Let $l$ be an ${\Bbb F}_q$-line in $H$.
Then ${}^{\#}(l \cap C ({\Bbb F}_q)) =0$ or $1$ or $d$ or $q+1$
by Proposition~\ref{possiblenumbers},
however ${}^{\#}(l \cap C ({\Bbb F}_q)) =q+1$ never occur because
$C$ has no ${\Bbb F}_q$-line components and $d <q+1$.

For each $i=0$ or $1$ or $d$,
let
\[
x_i = {}^{\#} \{l\,  :\mbox{\rm an ${\Bbb F}_q$-line $\subset H$}
   \mid  {}^{\#}(l \cap C ({\Bbb F}_q) ) =   i \}.
\]
Then
\[
x_0 + x_1 + x_d = \theta_q(2).
\]

Furthermore
${}^{\#}(l \cap C ({\Bbb F}_q)) =1$ if and only if
$l$ is the tangent line at an ${\Bbb F}_q$-point.
In fact, the ``if" part is obvious by the possible values of
${}^{\#}(l \cap C ({\Bbb F}_q))$. 
We see the ``only if" part of this claim.
Let $P_0$ be the ${\Bbb F}_q$-point of $l\cap C({\Bbb F}_q)$.
Suppose $l$ is not the tangent line $T_{P_0}(C)$
to $C$ at $P_0$.
Then
${}^{\#}(l \cap C ({\Bbb F}_q)) = {}^{\#}(T_{P_0}(C) \cap C ({\Bbb F}_q))=1.$
Let $m_1, \dots , m_{q-1}$ be the ${\Bbb F}_q$-lines passing through $P_0$
other than $l$ or $T_{P_0}(C)$.
Since
${}^{\#}(m_i \cap C ({\Bbb F}_q) \setminus \{ P_0 \}) \leq d-1$,
\[
N_q(C) \leq (d-1)(q-1) +1,
\]
which contradicts $N_q(C)= (d-1)q +1.$
Therefore $x_1=N_q(C) = (d-1)q +1$.

Consider the correspondence
\[
 \mathcal{P} = \{ 
 (P, l) \in C({\Bbb F}_q) \times \check{\Bbb P}^2({\Bbb F}_q)
 \mid {}^{\#}(l \cap C ({\Bbb F}_q) ) =   d, \ P \in l \},
\]
where $\check{\Bbb P}^2({\Bbb F}_q)$
is the set of ${\Bbb F}_q$-lines that are contained in
$H =  {\Bbb P}^2$.
Let $p_1: \mathcal{P} \to C({\Bbb F}_q)$
and $p_2: \mathcal{P} \to \check{\Bbb P}^2({\Bbb F}_q)$
be projections.
Since any ${\Bbb F}_q$-line $l$ on $H$ passing through
$P \in C({\Bbb F}_q)$ except the tangent line $T_P(C)$
meets $C({\Bbb F}_q)$ in $d$ points,
${}^{\#} p_1^{-1}(P) = q$.
So ${}^{\#}\mathcal{P} = ((d-1)q+1)q$.
On the other hand, ${}^{\#} p_2^{-1}(l) = d$
for any $l \in \image p_2$.
Hence,
\[
x_d = {}^{\#}\image p_2 = ((d-1)q+1)q/d.
\]

Therefore
\begin{align*}
 x_0 &= \theta_q(2) - x_1 - x_d \\
     &= -\frac{q}{d}(d -(\sqrt{q}+1))(d+\sqrt{q}-1),
\end{align*}
which must be nonnegative.
So we have $ d \leq \sqrt{q}+1$.
Combining this with the latter part of Lemma~\ref{mustbeirreducible},
we have $d=2$ or $\sqrt{q}+1$ if $d < q+1$.
\qed
\section{Uniqueness}
As we saw in the previous section,
the possible degrees of surfaces in Setting~\ref{oursetting}
are $2$, $\sqrt{q}+1$ and $q+1$.
\subsection{$d=2$}
This case is classical \cite{hir1985}.
Only hyperbolic quadric surfaces have $(q+1)^2$ points over ${\Bbb F}_q$,
and they are projectively equivalent to each other over ${\Bbb F}_q$.
\subsection{$d = \sqrt{q}+1$}
This case has been discussed in \cite{hom-kim2013preprint}.
\subsection{$d=q+1$}
The number $(d-1)q^2 + dq +1$ for $d=q+1$ is exactly $\theta_q(3)$.
This case was already handled by Tallini~\cite{tal}
in a little more general context.
\begin{proposition}
Let $S$ be a surface of degree $q+1$ in ${\Bbb P}^3$ over ${\Bbb F}_q$
without ${\Bbb F}_q$-plane components.
If $S({\Bbb F}_q) = {\Bbb P}^3({\Bbb F}_q)$,
then $S$ is projectively equivalent to the surface
\[
X_0X_1^q - X_0^qX_1 + X_2X_3^q - X_2^qX_3 =0
\]
over ${\Bbb F}_q$.
\end{proposition}
\proof
In ${\Bbb P}^3$ with homogeneous coordinates $X_0, \dots , X_3$,
the ideal of the algebraic set ${\Bbb P}^3({\Bbb F}_q)$
is generated by
\[
\{ X_iX_j^q - X_i^qX_j \mid 0 \leq i < j \leq 3 \}.
\]
Hence $S$ is defined by an equation of the type
\[
(X_0, \dots , X_3) A \, {}^t(X_0^q, \dots , X_3^q)= 0,
\]
where $A$ is an alternating matrix over ${\Bbb F}_q$.
When $q$ is a power of $2$,
we understand an alternating matrix over ${\Bbb F}_q$
to be a symmetric matrix whose diagonals are $0$.
By the standard argument of linear algebra,
we can find a nonsingular matrix $G$ over ${\Bbb F}_q$
such that
${}^tGAG$ is either
\[
\begin{pmatrix}
 \begin{matrix}
    0 & 1 \\
    -1& 0
 \end{matrix} & {\rm \Large O} \\
   {\rm \Large O} & {\rm \Large O}
\end{pmatrix}
\ 
\mbox{\rm or}
\ 
\begin{pmatrix}
 \begin{matrix}
    0 & 1 \\
    -1& 0
 \end{matrix} & {\rm \Large O} \\
   {\rm \Large O} & \begin{matrix}
                      0 & 1 \\
                      -1& 0
                    \end{matrix}
\end{pmatrix}.
\]
Note that since all entries of $G =(g_{ij})$ are in ${\Bbb F}_q$,
the matrix $G^{(q)} =(g_{ij}^q)$ coincides with $G$.
Therefore after changing coordinates by $G$,
we get an equation of $S$ as
\[
X_0X_1^q - X_0^qX_1=0
\]
or
\[
X_0X_1^q - X_0^qX_1 + X_2X_3^q - X_2^qX_3 =0,
\]
but the former splits into $q+1$ ${\Bbb F}_q$-linear polynomials.
Hence only the latter case occurs.
\qed

\gyokan

\noindent
{\bf Acknowledgment}
This work was done during our stay at
Minnesota state University, Mankato.
We thank Professor Namyong Lee for his hospitality.


\end{document}